\documentclass[a4paper, 11pt]{article}

\usepackage{import} 
\usepackage[left=2.5cm,right=2.5cm,top=2.5cm,bottom=2.5cm]{geometry}

\usepackage{subfiles}

\usepackage[authoryear, square, sectionbib]{natbib}

\usepackage{amssymb, amsthm, mathtools}
\usepackage{dsfont} 
\usepackage[ruled,linesnumbered]{algorithm2e}

\usepackage{booktabs} 
\usepackage{multirow} 
\usepackage{diagbox} 

\usepackage{graphicx} 
\usepackage[dvipsnames]{xcolor} 

\usepackage{caption} 
\usepackage{subcaption}
\usepackage{wrapfig} 

\usepackage{silence} 

\usepackage{hyperref} 
\usepackage[noabbrev]{cleveref}

\newcommand{\myBullet}{$\bullet~$}

\AtBeginEnvironment{subappendices}{%
    \section*{Appendix}
    \addcontentsline{toc}{section}{Appendices}
}

\newcommand{\mc}{\mathcal} 

\newcommand{\simiid}{\overset{i.i.d}{\sim}}
\newcommand{\siminde}{\overset{i.}{\sim}}

\newcommand{\som}[3]{\sum_{#1=#2} ^{#3} }
\newcommand{\pro}[3]{\prod_{#1=#2} ^{#3} }

\newcommand{\pen}{\mathrm{pen}}
\newcommand{\regularization}{\lambda}

\DeclareMathOperator{\prox}{Prox}
\DeclareMathOperator*{\argmin}{arg\,min}
\DeclareMathOperator*{\argmax}{arg\,max}

\newcommand{\setr}{\mathbb R}

\newcommand{\sets}{\mathbb S}
\newcommand{\setsppr}[1]{\sets^{++}_{#1}}

\newcommand{\abs}[1]{\left\vert #1 \right\vert}
\newcommand{\norm}[1]{\left\|#1\right\|}

\newcommand{\sca}[1]{\left\langle#1\right\rangle}

\newcommand{\p}[1]{{\left(#1\right)}}

\newcommand{\stackeq}[2][rrrrrrrrr]{\left\{\begin{array}{#1} #2 \end{array}\right.}

\newcommand{\memfct}{m}
\newcommand{\marg}{g}
\newcommand{\comp}{f}

\newcommand{\para}{\theta} 
\newcommand{\paraspace}{\Theta} 

\newcommand{\nobs}{N} 
\newcommand{\nrep}{J} 

\newcommand{\obs}{Y} 

\newcommand{\lat}{Z} 

\newcommand{\varlat}{\Omega} 
\newcommand{\meanlat}{\mu} 

\newcommand{\covrep}{t} 
\newcommand{\eps}{\varepsilon} 
\newcommand{\vareps}{\Sigma} 
\newcommand{\dimparind}{q} 

\newcommand{\covgd}{X} 
\newcommand{\dimcovgd}{p} 
\newcommand{\dimld}{a} 

\newcommand{\iter}{k}

\newcommand{\step}{\gamma}

\newcommand{\Prec}{P}
\newcommand{\invPrec}{\Prec^{-1}}

\theoremstyle{plain}
\newtheorem{theorem}{Theorem}[section]

\theoremstyle{definition}
\newtheorem{definition}[theorem]{Definition}

\theoremstyle{remark}
\newtheorem{remark}[theorem]{Remark}

\newtheoremstyle{exampstyle}
    {\topsep}
    {\topsep}
    {\itshape}
    {}
    {\bfseries}
    {.}
    {.5em}
    {}

\theoremstyle{exampstyle}

\theoremstyle{remark}

\SetKwInOut{Require}{Require}
\SetKw{Initialize}{Initialize}
\SetKw{Return}{Return}

\SetKwFor{Tq}{tant que}{faire}{fin}
\SetKwInput{Entree}{Entrées}
\SetKw{Retour}{retourner}

\newcommand{\codename}[1]{\textsc{#1}\xspace}

\newcommand{\R}{\codename{R}}

\newcommand{\saemix}{\codename{Saemix}}
\newcommand{\nlm}{\codename{nlm}}
\newcommand{\glmnet}{\codename{glmnet}}
\newcommand{\glmmLasso}{\codename{glmmLasso}}

\newcommand{\PSPG}{\codename{AWPSG}}

\newcommand{\NLMEMLASSO}{\codename{NLMEM-LASSO}}

\newenvironment{figureht}{
    \begin{figure}[!ht]\centering
	}{\end{figure}}
		
\newenvironment{tableht}{
    \begin{table}[!ht] \centering
	}{\end{table}}

\author{Antoine Caillebotte}
\title{Estimation and variable selection in high dimension  in nonlinear mixed-effects models.}

\hypersetup{
    pdfauthor={Antoine Caillebotte, Estelle Kuhn, Sarah Lemler},
    pdftitle={Estimation and variable selection in high dimension  in nonlinear mixed-effects models.},
    pdfkeywords={Statistic, High dimension, Variable selection}
} 

\newcommand{\mailto}[1]{ \href{mailto:#1}{#1}}

\begin{document}

    \begin{center}
        {\Large
        	{\sc  Estimation and variable selection  in high dimension in nonlinear mixed-effects models.}
        }
        \bigskip
        
       Antoine Caillebotte $^{1,2}$ \& Estelle Kuhn  $^2$ \& Sarah Lemler $^3$   
        \bigskip
        
        {\it
            $^1$ Université Paris-Saclay, INRAE, UMR GQE-Moulon, France,      \mailto{caillebotte.antoine@inrae.fr}, \\
            $^2$ Université Paris-Saclay, INRAE, UR MaIAGE, France,
            \mailto{estelle.kuhn@inrae.fr},\\
            $^3$ Université Paris-Saclay, CentraleSupélec, Laboratoire MICS, France,
            \mailto{sarah.lemler@centralesupelec.fr}
        }
    \end{center}
\bigskip
{\bf Abstract.}
    We consider nonlinear mixed effects models including high-dimensional covariates to model individual parameters variability. The objective is to identify relevant covariates among a large set under sparsity assumption and to estimate model parameters. To face the high dimensional setting we consider a regularized estimator namely the maximum likelihood estimator  penalized with the $\ell_1$-penalty. We rely on the use of the eBIC model choice criteria to select an optimal reduced model. Then we estimate the parameters by maximizing the  likelihood of the reduced model.  We calculate in practice the maximum likelihood estimator  penalized with the $\ell_1$-penalty though a weighted proximal stochastic gradient descent algorithm with an  adaptive learning rate. This choice allows us  to consider very general models, in particular models that do not  belong to the curved exponential family. We demonstrate first in a simple linear toy model through a simulation study the good convergence properties of this optimization algorithm. We compare then the performance of the proposed methodology with those of the \glmmLasso procedure in a linear mixed effects model in a simulation study. We  illustrate also  its performance  in a nonlinear mixed-effects logistic growth model through simulation.  We finally highlight the beneficit of the proposed procedure relying on  an integrated single step approach regarding two others two steps approaches for variable selection objective.

{\bf Keywords.} nonlinear mixed effects model,
     high dimension, variable selection, LASSO penalty, stochastic gradient descent, weighted proximal.

\newpage

    \pagenumbering{arabic}

\section{Introduction}

    Mixed effects models are very fine and very useful statistical modeling tools for analyzing data with hierarchical structures and repeated measurements (see \cite{pinheiro_mixed_effects_2006}). In particular it is possible to account for several levels of variability of a phenomenon observed within a population of individuals,  each individual being observed several times in different conditions.  They are widely used in many applied fields such as agronomy, pharmacology, and even economics. Mixed effects models are composed of two nested levels of modeling: on the one hand a common structural modeling of the phenomenon of interest parameterized for each individual in the population by specific individual parameters, on the other hand a modeling of these individual parameters as random variables which accounts for their variability within the population. The first modeled the intra individual variability among the repeated measurements of each individual, the latter stands for the inter individual variability between individuals in the population.
    
    The structural modeling of the phenomenon of interest can be linear or non-linear in the individual parameters. Non-linear type modeling can in particular account for complex phenomena modeled mechanistically by models which integrate, for example, physical or biological knowledge. The parameters of such models are often of strong practical interest because they are interpretable from an applied point of view. In such setting a central objective is to characterize and explain the variation of these individual parameters within the population. Therefore the modeling of individual parameters is done through random effects at the individual level and can also integrate descriptive covariates of the individuals.
    
    Depending on the context, these descriptive covariates can be of high dimension and the objective is then to identify those which are the most relevant to explain the variabilities observed within the population by estimating the associated vector of regression parameters. Let us consider the field of plant ecophysiology. Mechanistic models have been proposed to describe plant development processes. These models integrate descriptive variables of the environment as covariates acting directly on the plant development process. The parameters of these models are often physical quantities such as leaf appearance speeds or light interception capacities. These parameters may vary when considering a population of plants from different varieties each characterized by its genotype (see \cite{baey2018mixed}). Their variations can be modeled into a mixed effects model and also integrate large genetic markers characterizing the genotypic variability within the population. In this context, identifying the relevant covariates among a set of high-dimensional covariates amounts to identifying the genetic markers that influence the phenomenon of interest (see \cite{Onogi2020}).
    
    From the point of view of inference in mixed effects models including high-dimensional covariates, the objective is on the one hand to select the relevant covariates from a set of high-dimensional covariates in order to identify a parsimonious model with a reduced number of covariates and on the other hand to estimate the parameters in the reduced model. There are two main difficulties for the inference task. The first one is the high dimension of the covariates.  Usually the  selection of relevant covariates can be done via a regularization approach. Furthermore, in the context of mixed effects models, an additional difficulty appears due to the presence of random effects which are not observed. This is a classic context of latent variable models. Inference is complex to carry out in this framework due to the latent structure and involves the use of efficient numerical methods. Parameter inference can be done for example by maximum likelihood via Expectation Maximization (EM)  or stochastic gradient type algorithms. In the context of exponential family models, EM type algorithms are easy to implement and have good theoretical properties (see \cite{dempster1977maximum, delyon1999convergence}). On the other hand, to our knowledge, there are no theoretical results of convergence outside this framework. In addition, the implementation of these algorithms outside the exponential family is more complex. To get around this limitation, a trick called exponentialization trick is sometimes used in practice. However, its limits have been highlighted in \cite{debavelaere2021}. In particular, this procedure can generate significant estimation biases due to the fact that the inference of the parameters is carried out in an extended model different from the initial one. Stochastic gradient methods can be applied in more generic models, in particular outside the exponential family. Theoretical guarantees of convergence towards an extremum of the target function have been established under general assumptions.  However, these methods are quite  sensible   in practice to the tuning of the learning rate,  in particular in high-dimensional parameter spaces due to the heterogeneity of the different components of the gradient. Adaptive choices of the learning rate  and procedures based on gradient preconditioning have been proposed in generic contexts to overcome this computational difficulty \citep{Duchi_AdaGrad_2011, Diederik_Adam_2017}. More recently in the context of maximum likelihood estimation in general latent variable models, a stochastic gradient algorithm integrating a gradient preconditioning step based on an estimator of the Fisher information matrix obtained as a by-product  of the algorithm has  been proposed, opening new possibilities for maximum likelihood inference (\cite{SGDFIM_Baey_2023}).
    
    In order to achieve the objective of selecting relevant variables from a set of high-dimensional variables in mixed effects models, several regularization approaches have been developed. \cite{GLMMLasso_Schelldorfer_2011}  proposed a maximum likelihood estimator with a LASSO penalty (\cite{tibshirani_regression_1996}) in the context of linear mixed effects models with a Gaussian error term and developed a R package (see \cite{schelldorfer2014glmmlasso}). \cite{fort2019stochastic} and  \cite{ollier2022fast} proposed estimators with more general penalties in the case of nonlinear mixed effects models belonging to the curved exponential family. \cite{bertrand2013multiple} have also proposed an stochastic penalized versions of the EM (\cite{delyon1999convergence}) algorithm to take into account multiple parameters in pharmacokinetic models. Bayesian approaches based on “spike and slab” distributions have also been developed by \cite{HeuclinBayesianSelection_2020} in the case of  linear mixed models and by \cite{naveau2024bayesian} in the non linear case. On the other hand, to our knowledge, there are no high-dimensional variable selection methods for general nonlinear mixed-effects models,  in particular  outside the exponential family.
    
    In this contribution, we propose a variable selection and re-estimation procedure in a general mixed effects model with high dimensional covariates explaining individual parameters variability. We deal with the variable selection in high dimension through a maximum likelihood estimate regularized via a $\ell_1$-penalty term.  To calculate this penalized estimator, we use an adaptive  weighted proximal stochastic gradient algorithm to simultaneously handle the presence of latent variables and the non differentiability of the penalty term.  The latter allows to handle very general models, in particular   our procedure does not require that the  models   belong to the curved exponential family. This paper is organised as follows. The second section introduces the mixed-effects models with high-dimensional covariates and presents two classical examples. The third section is devoted to the description of the proposed  variable selection  procedure in  the high-dimensional setting and the parameter's re-estimation in the reduced model. The fourth section gathers practical details of the  numerical methodology.  Finally, we present a simulation study  and discuss the potential of the proposed method with respect of existing ones. The paper ends with some conclusion and perspectives.

\section{Nonlinear mixed effects model integrating high dimensional covariates to model individual variability}

    \subsection{Model description}
        Let $\nobs$ be a positive integer. We consider  $\nrep$ repeated measurements for each individual $i$. Therefore we have  $\nrep$ observations per individual. Let us  denoted by $\obs_{i,j}$ the $j-$th observation of the $i-$th individual for $1\leq i\leq\nobs$ and $1\leq j\leq \nrep$. We assume that $\obs_{i,j}$ takes value in $\setr^d$. We model this  observation  with a nonlinear mixed effects model (see \cite{pinheiro_mixed_effects_2006} and \cite{davidian2017nonlinear}):        
        \begin{align*}
            \left\{
                \begin{array}{ll}
                    \obs_{i,j}&=\memfct(\alpha,\covrep_{i,j},\lat_i)+\eps_{i,j}, \\
                    \lat_i&\siminde\mathcal{N}(\covgd_i\beta,\varlat), \\
                    \varepsilon_{i,j}&\simiid\mathcal{N}(0,\vareps).
                \end{array}
            \right.
        \end{align*}
        where $\memfct$ is a  function  depending on population parameter  $\alpha$ taking values in $\setr^a$, individual covariates $(\covrep_{i,j})$ and individual parameters modeled by the latent variable $\lat_i$. We consider functions that can be either  linear or nonlinear with respect to the latent variable. 
        Note that in the case of longitudinal data, the covariates $\covrep_{i,j}$ stands usually  for the  $j-$th observation time of individual $i$.  Individual parameters $(\lat_i)$    describes the inter-individual variability of the population. 
        For the $i-$th subject, the quantity $\lat_i$ is a   $\dimparind-$dimensional random vector, independent of $(\varepsilon_{i,j})_j$ and assumed to be distributed as a Gaussian distribution with expectation $\covgd_i\beta$ and covariance matrix $\varlat$, where $\covgd_i$ is a matrix of size $\dimparind \times \dimcovgd$ and $\beta$ is a vector of size $\dimcovgd$. The noise term $\varepsilon_{i,j}$ is usually assumed centered Gaussian with unknown covariance matrix $\vareps$. The unknown parameters of the nonlinear mixed-effects model are therefore $\para=(\alpha,\beta, \varlat,\vareps)\in\setr^\dimld\times \setr^{\dimcovgd}\times \setsppr{\dimparind}\times\setsppr{d}$ where $\setsppr{d}$ stands for the set of symmetric positive definite matrices of size $d \times d$.
        We emphasize that our model formulation is very general and allows us to consider model function $\memfct$ depending on both population parameters and individual parameters. We consider the setting of high-dimensional covariates where $\dimcovgd$ can be much larger than $\nobs$ and assume that only a small subset of covariates are relevant to be included in the model.
        
        \begin{remark}
            More general distributions than the Gaussian one can be chosen for the noise term as well as more general settings with various numbers of observations per individual can  be considered.
        \end{remark}
        
        \begin{remark}
            It may be useful in practice to reparametrize the latent variables $(\lat_i)$ and/or the parameters. For example, consider as individual parameter $\log(\lat_i)$ rather than $\lat_i$ if individual parameters are assumed to be positive. We refer for more details to the  reparametrization cookbook (\cite{paramcookbook}).
        \end{remark}

    \subsection{Practical examples}

We present in this section two classical nonlinear mixed effects models, one modeling general growing process and one from the pharmacological field.
    \subsubsection{Logistic growth curves}
    
        We consider as  first example the specific case of the logistic  curve model used to model growth dynamic, which is commonly used in the nonlinear mixed-effect models' community and presented in \cite{pinheiro_mixed_effects_2006}. The observation $Y_{i,j}$ is the circumference of the tree $i$ at time $t_{ij}$ and is modeled througth the logistic model
     with the function $\memfct$  given by:
        \begin{equation}\label{eq:mem_orangetrees}
            \memfct(\alpha,\covrep_{ij},\lat_i) = \frac{\lat_{i1}}{1 + \exp\p{-\frac{\covrep_{ij} - \lat_{i2}}{\alpha}}},
        \end{equation}
        where the individual parameters  $\lat_{1i}$ represents the asymptotical maximum value of the circumference, $\lat_{2i}$ represents the value of the sigmoid's midpoint, and $\alpha$ represents the logistic growth rate. 

    \subsubsection{Pharmacodynamic model}
        
        We consider as a second  example the two compartments pharmacodynamic model used by \cite{pinheiro_mixed_effects_2006}. The observation $\obs_{i,j}$ is the serum concentration measured at time $t_{i,j}$. The model is given by:
        \begin{equation}\label{eq:mem_theophylline}
            \memfct(\covrep_{ij},\lat_i) = \frac{D_i k_{ai}}{V_i(k_{ai} - Cl_i/V_i)} \left(\exp(-k_{ai}\covrep_{ij}) - \exp\left(-\frac{Cl_i}{V_i}\covrep_{ij}\right)\right),
        \end{equation}
        where the individuals parameters $V_i$ represents the distribution volume, $k_{ai}$ the absorption rate,  $Cl_i$  the clearance and $Cl_i/V_i$ the elimination rate; $D_i$ is the known dose of the drug receive by individual  $i$.

\section{Variable selection and estimation procedure}
    
    In this section, we present the proposed method called \NLMEMLASSO for both  selecting the relevant covariates and estimating the model parameters.

    \subsection{Estimation in latent variable model }
            
        We consider the maximum likelihood estimator to infer the non-linear mixed-effects model's parameters. In the context of latent variable models, the marginal likelihood, denoted by $\marg$, is obtained by integrating the complete likelihood over the latent variables, which are not observed.
        \begin{eqnarray}
            \marg(\theta;\obs) 
            &=& \pro i1{\nobs}\int f(\theta;\obs_i,\lat_i) d\lat_i \notag
             =\pro i1{\nobs}\int  p_\theta(\obs_i|\lat_i)p_\theta(\lat_i) d\lat_i \notag \\
            &=& \pro i1{\nobs}\int \left\{ \pro j1\nrep p_\theta(\obs_{i,j}|\lat_i)\right\}p_\theta(\lat_i)d\lat_i \label{eq:Lcomp}
        \end{eqnarray}
        
        \noindent where $\comp(\theta;\obs,\lat)$, $p_\theta(\obs|\lat)$, $p_\theta(\lat)$ are respectively  the density of the pair $(\obs,\lat)$, the density of $\obs$ conditionally to $\lat$, and the density of $\lat$.
               
        One can usually  estimate the model parameters by maximizing the logarithm of the marginal likelihood using the maximum likelihood estimator  written as follows:
        
        \begin{equation}
            \label{eq:MLE_estimator}
            \hat\theta^{\text{MLE}} = \argmax_{\theta\in \paraspace} \log \marg(\theta;\obs) ,
                \end{equation}
        where $\paraspace$ denotes the parameter space. From a practical point of view, this estimate can often not be calculate explicitly. To deal with latent variables, classical methods used to infer the unknown parameters are Expectation Maximization like algorithms (see \cite{ng2012algorithm}). However two major limits of these procedures is that i) theoretical convergence properties are established only for models belonging to the curved exponential family; ii)
               their implementation is easy for   models of the curved exponential family but much tricky outside this setting.  However one can use the exponentiation trick by defining certain parameters as Gaussian variables and study an augmented model that belongs to the exponential family.  Nevertheless, its limits have been highlighted in \cite{debavelaere2021}.  Another type of numerical approaches to calculate the MLE are stochastic gradient descent algorithms (see \cite{cappe2005springer}).
                    Recently \cite{SGDFIM_Baey_2023} have presented a preconditioned stochastic gradient descent for estimation in a latent variable model adapted to general latent variables models, in particular theoretical  guarantees are obtained for general models, without assuming that the models belong to the curved exponential family.
    
    \subsection{Penalized likelihood in high-dimensional setting}
        
        In our context, we must deal with the high dimension of the covariates, therefore we introduce a regularization term  and consider a penalized maximum likelihood estimator. We aim to select relevant variables among the covariates and follow the idea of  the LASSO (Least Absolute Shrinkage and Selection Operator) procedure, which was initially developed for linear regression models in  \cite{tibshirani_regression_1996}. This method enables us to handle high-dimensional covariates and select a subset of explanatory covariates from a large collection. Therefore we introduce  a $\ell_1$  penalty term, which only depends on the parameter $\beta$:                
        \begin{equation}\label{eq:LASSO}
            \pen_\regularization(\theta) = \regularization \norm{\beta}_1 = \regularization \som k1p \abs{\beta_k},
        \end{equation}
        where $\regularization$ is a positive real called  regularization parameter. Our goal is then to maximize
        the penalized criterion defined as the difference  of the logarithm of the marginal likelihood and of  the penalty term. Let us define the penalized maximum likelihood estimator by:        
        \begin{equation}\label{eq:penalized_estimator}
            \hat\theta^{PEN}_\regularization= \argmax_{\theta\in \paraspace} \left\{\log \marg(\theta;\obs) - \pen_\regularization(\theta)\right\},
        \end{equation}
        where $\paraspace$ denotes the parameter space and where $\regularization$ is a positive parameter. The larger the value of $\regularization$, the more $\beta$ will be constrained to have zero components. Conversely, the smaller the value of $\regularization$, the freer the components of $\beta$ will be. It is customary to determine the value of $\regularization$ using cross-validation  or using model criterion (see \cite{tibshirani_regression_1996}). We will consider the eBIC model choice criterion well adapted to the high dimensional setting to find on optimal  regularization value (see \cite{eBIC_chen_chen_2008}).

    \subsection{Regularization path procedure}
    
        To select only the most explanatory covariates, it's important to choose a well-balanced value for the regularization parameter $\regularization$. We choose an optimal parameter by minimizing the extended Bayesian Information Criterion (eBIC) (see \cite{eBIC_chen_chen_2008}).  Guided by the intuition given by \cite{delattre_mem_BIC_2014}, we penalize the number of degrees of freedom by the log of the total number of observations $\nobs \nrep$. 
        We consider  a grid  $\Lambda$ of values for the regularization parameter $\regularization$. We conduct then the following inference methodology \label{algo:metho}: 
            
        \begin{enumerate}
            \item[i)] For all $\regularization\in\Lambda$ repeat the following steps:
            \begin{itemize}
\item            calculate
                $
                \hat\theta^{PEN}_\regularization= \argmax_{\theta\in \paraspace} \left\{\log \marg(\theta;\obs) - \pen_\regularization(\theta)\right\}.
                $
                \item  deduce the associated support $\hat S_\regularization=\{j\in\{1,\dots,p\}, \hat\beta_{\regularization,j}^{PEN}\neq 0\}$
                \item calculate  $
                \hat\theta^{\text{MLE}}_\regularization= \argmax_{\theta\in \paraspace} \left\{\log \marg_\regularization(\theta;\obs) \right\},
                $ where  $\marg_\regularization$ is the marginal likelihood in the  model restricted to the support  $ \hat S_\regularization$.
                \item calculate the eBIC criterion:
                \begin{equation} \label{eq:eBIC}
                    \text{eBIC}(\regularization)=-2\log \marg_\regularization(\hat\theta^{MLE}_\regularization; \obs)+|\hat S_\regularization|\log(\nobs\nrep) +  2\log\left(\binom{\dimcovgd}{|\hat S_\regularization|}\right)
                \end{equation}
                \end{itemize}

 \item[ii)] Select the parameter $\regularization$ that minimizes the eBIC:
                \[\hat\regularization=\argmin_{\regularization\in\Lambda}\text{eBIC}(\regularization),
                \]
                and consider the final estimator  defined by $\hat\theta^{MLE}_{\hat\regularization}$.
        
        \end{enumerate}

\section{Numerical methodology in practice}

In order to carry out in practice our variable selection and re-estimation procedure, we need to be able to efficiently and accurately calculate both the maximum likelihood estimate in a reduced model and the penalized estimate for numerous different values of the regularization parameter. For both objectives, we have to deal with the integral on the latent variables, which means that there are no explicit expression for the criteria to be optimised in most cases.
 Therefore, we will consider numerical methods to solve these maximization problems. We propose to use an adaptive stochastic gradient algorithm for the first objective and a proximal weighted extended version for the second objective.

    \subsection{Adaptive Stochastic Gradient Descent Algorithm}

            Adaptive algorithms such as AdaGrad \citep{Duchi_AdaGrad_2011}, RMSPprop \citep{Tieleman_RMSP_2012} and Adam \citep{Diederik_Adam_2017} have proved their worth in many settings. The objective of these strategies is to automatically rescale through the adaptive learning rate each gradient descent direction, in order to scale the different components of the gradient, homogenizing the evolution of the algorithm. We choose to benefit from the advantages offered by the adaptive algorithm AdaGrad. The main advantage is the adaptive property of the proposed learning rate, since the convergence behavior of stochastic gradient algorithms depends strongly of the choice of the learning rate. The second advantage relies on the component wise independence property of the Adagrad learning rate which will be of great importance  when  dealing with the non-differentiable penalty term of the criterion in the next section.   
            The algorithm is divided into two steps. The first step consists in the simulation of  a realization of the latent variables, through a direct sampling from the posterior distribution using the current value of the parameter or a step of a Metropolis-Hastings sampler (see \cite{marin2007bayesian}). The second step is the update of the parameter estimate using the adaptive preconditioned gradient descent on the approximate complete likelihood, More precisely, given the realization of the latent variable, we evaluate the gradient of the complete log-likelihood and calculate for each component of the parameter the cumulative sums of the square of these gradients, and define the preconditioning diagonal matrix proportional to the inverse of  the square root of these values.

           \subsection{Adaptive  Weighted Proximal Stochastic Gradient Descent Algorithm} 
\label{sec:proximalweighted}

             When evaluating the penalized estimate in practice,  we have to face to the 
        non-differentiability of the $\ell_1$-penalty term and consider therefore a proximal algorithm  as presented for example by \cite{fort2019stochastic}. 
            Note that in our context of high-dimensional setting, it is difficult to manage  a full matrix as preconditioner from a computational point of view. Therefore the adaptive vectorial learning rate detailed in the previous section is particularly well-adapted since the preconditioning matrix is diagonal and thus allows to evaluate the proximal step component by component.           
           The main idea is to combine the stochastic gradient descent algorithm  with a Proximal Forward-Backward algorithm (\cite{Chen_ConvergeFB_1997, Tseng_FB_2000}).   The algorithm is divided into three steps;  a realization of the latent variables is sampled with a first step called \textit{Simulation}, as in the previous case. The second step is the adaptive preconditioned gradient descent on the approximate complete likelihood, called the \textit{Forward} step. The last step, called \textit{Backward}, deals with the penalty term. We apply the proximal operator (\cite{moreau_fonctions_1962,rockafellar_monotone_1976, VMFB_Chouzenoux_2014}) detailed in the following paragraf.

           \renewcommand{\Prec}{A}
       
        We first describe the weighted proximal operator.
            First of all, for any symmetric positive-definite matrix $\Prec\in\setsppr{d}$, let's denote by $\norm ._\Prec$ the $A$-weighted norm defined with the $\Prec$-weighted scalar product by $\forall (x,y)\in\setr^d, \sca{x,y}_\Prec = \sca{x,\Prec x}$ and $\norm{x}^2_\Prec = \sca{x,x}_\Prec$.
    
            \begin{definition} Assume that $\Prec \in\setsppr{d}$. Let $\psi : \paraspace \rightarrow \setr$ be a proper, lower semicontinuous, convex function and $\theta \in\paraspace$. The weighted proximal operator of $\psi$ in $\theta$ relative to the metric induced by $A$ is the unique minimizer of $\psi(\theta') +\frac 12 \norm{. - \theta}_\Prec^2$:
                \begin{eqnarray}
                    \prox_{\Prec, \psi}(\theta) = \argmin_{\theta'\in\setr^d}\left( \psi(\theta') +\frac 12 \norm{\theta' - \theta}_A^2\right)
                \end{eqnarray}
            \end{definition}

            The  general resulting update step of the proximal gradient descent algorithm targeting the minimum of a given function $f$ can then be written as:
            \begin{equation*}
                \theta_{\iter+1} = 
                \prox_{\Prec_\iter,\pen_{\regularization}}(\theta_{\iter}-\invPrec_\iter\nabla_\theta\comp(\theta_\iter))
            \end{equation*}
            where $(\Prec_\iter)_{\iter\geq 1}$ is a sequence of symmetric positive-definite preconditioning matrices.

 \begin{remark}
                We emphasize that when  $\Prec_\iter = \gamma^{-1}_\iter Id$, the operator $\prox_{\gamma^{-1}Id, \psi} $ is the classical isotropic proximal operator defined by \cite{moreau_fonctions_1962} and we obtain the standard proximal stochastic gradient descent where $(\gamma_k)_{k\geq 1}$ is the common step size. 
            \end{remark}

            In general, $\prox_{\Prec,\pen}$ is well defined but may not be quickly computable. \citeauthor{Becker_prox_computable_2012} state several assumptions on $\Prec$ and the function $\pen$ that allow the proximal operator to be computable as for  example, rank restriction and $\pen$ is separable, i.e. $\pen(x) = \som l 1\dimcovgd \pen_l(x_l)$. If no assumptions are made, the proximal operator is often not explicit and will require additional steps. Since we choose the AdaGrad type learning rate and the $\ell_1$-penalty, the calculation of the proximal operator is separable and can be perfomed component per component independently and has  an explicit form: for a given $\regularization>0$, we get for all $s>0$ and $x \in \mathbb{R}$ \renewcommand{\Prec}{s}
            \begin{equation}\label{eq:prox_lasso}
                \prox_{\Prec, \regularization|.|}(x) = \stackeq[cr]{
                    0 &\text{if } \abs{x} < \regularization / \Prec
                \\ x - \regularization / \Prec
                    &\text{if  } x \geq \regularization / \Prec
                \\ x + \regularization / \Prec
                    &\text{if } x \leq -\regularization / \Prec
                    } 
            \end{equation}
We provide  here in detail the steps of the adaptive weighted proximal stochastic gradient descent algorithm.

            \begin{algorithm}[!ht] \label{algo:SPG-fim}
                \caption{Adaptive  Weighted Proximal Stochastic Gradient (\PSPG)}
                \Require{$\step_0>0$, $\epsilon>0$, $K_{max} \in \mathbb{N}$ } 
                \Initialize Initialize starting point $\theta_0 \in \mathbb{R}^d$; adaptive stepsize vector $\Prec_0 =0$; iteration $k=0$.
                
                \While{$(\theta_\iter)$ not converged or $k \leq K_{max} $}{
                    $\iter \leftarrow \iter +1$
                    \\
                    \myBullet {\bf Simulation step :} \\ 
                    \quad Draw $\lat^{(\iter)}$
                    from the posterior distribution $p( . |  \theta_k)$ or
                    using a single step of a Hastings Metropolis procedure having as stationary distribution the  posterior.
                    \\
                    \myBullet {\bf Gradient computation :}
                    $v_\iter \leftarrow \frac 1\nobs \som i1\nobs \nabla \log \comp_{\theta_\iter}(\obs_i,\lat^{(\iter)}_i)$
                    \\
                    \myBullet {\bf Adaptive stepsize :}
                    $\Prec_{\iter,l}^2 \leftarrow \Prec_{\iter-1,l}^2+ v_{\iter,l}^2$ for $l \in\{1,\dots, \dimcovgd\}$
                    \\
                    \myBullet {\bf Update parameters :} \\
                    \quad \myBullet {\bf Forward step :} 
                    $\omega_{\iter,l} \leftarrow \theta_{\iter-1,l} + \step_0 v_{\iter,l}/(\Prec_{\iter,l}+\epsilon)$ for $l \in\{1,\dots, \dimcovgd\}$
                    \\
                    \quad \myBullet {\bf Backward step :} 
                    $\theta_{\iter,l} \leftarrow \prox_{(\Prec_{\iter,l}+\epsilon)/\step_0, \regularization|.|}(\omega_{\iter,l})$     for $l \in\{1,\dots, \dimcovgd\}$
                    
                }
                \Return{$\hat \theta = \theta_\iter$}
            \end{algorithm}

        \begin{remark}
            Regarding the simulation of the realization of the latent variable $\lat^{(\iter)}$ at iteration $\iter$, it is often the case that direct sampling from the posterior distribution is not possible and one must require to Monte Carlo Markov Chain procedure such as the Hastings Metropolis algorithm or the Gibbs algorithm (see \cite{marin2007bayesian}). 
        \end{remark}
        
        \subsection{Details on the regularization path}
\begin{figure}[!ht]
            \centerline{\includegraphics[width = 0.9\textwidth]{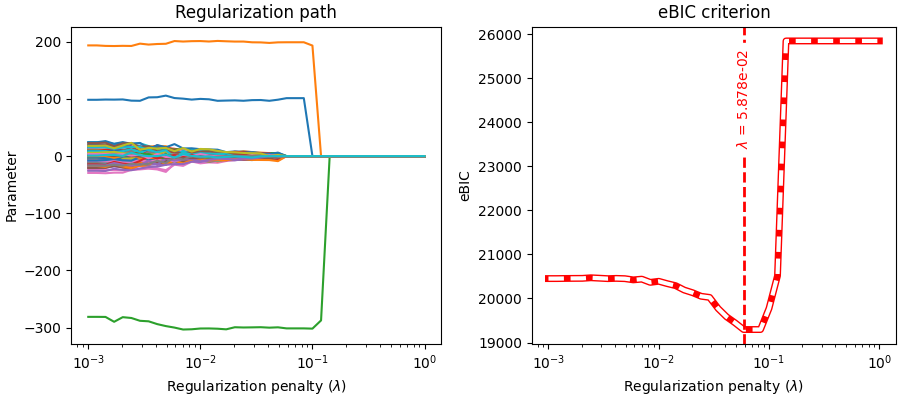}}
            \caption{Left: Regularization path, i.e. values of estimates of $\beta$, in solid line; Right: values of the eBIC criterion in dotted line; the dotted vertical line represents the chosen regularization value.}
            \label{fig:regularization_path}
        \end{figure}

            In this section, we give more detail about the procedure used to choose a regularization parameter $\regularization$ (see \cref{algo:metho}). The minimization of the eBIC criterion (\ref{eq:eBIC})  calculated on a  proposal grid $\Lambda$ can be represented by the  regularization path which can be plotted to visualize the evolution of the support of $\beta$ with respect to the regularization parameter $\regularization$ (e.g. \cref{fig:regularization_path}). This type of graph shows that the larger the regularization parameter, the more the component of the vector $\beta$ is restricted to zero. The smaller $\regularization$ is, the more $\beta$ components are free. In practice, the  proposal grid $\Lambda$ is built on a logarithmic scale. In fact, it's preferable to explore the smaller values a little more, in order to capture rapid changes in the support of $\beta$. This way, we don't spend too much time on ranges of regularization values that would have returned the same support.

        Note also  that to evaluate the eBIC criterion defined in (\ref{eq:eBIC}) one must evaluate the quantity $\marg_\regularization(\hat\theta^{MLE}_\regularization; \obs)$ which does not have an explicit form due to the presence of the latent variables $\lat$ and the non-linearity of the model. Therefore we use a Monte Carlo procedure to calculate an approximation of $\marg_\regularization(\hat\theta^{MLE}_\regularization; \obs)$.

\section{Numerical experiments}

    In this section, we study the performance of the proposed procedure through  simulations.  The numerical study is divided into three parts. First, we demonstrate  the good behavior of the optimization algorithm presented in section \ref{sec:proximalweighted}. Then we compare the performance  of the variable selection procedure with re-estimation in a linear mixed-effects model to those of the existing method \glmmLasso. We assess the robustness of our method when dealing with a nonlinear mixed-effect model in different configurations. We emphasize that to the best of our knowledge, there exists actually no other methodology to deal with this issue in general nonlinear mixed-effect model.  Finally, we highlight the benefits of the proposed integrated procedure, which combines the mixed population model and the variable selection procedure in a single step, compared with two other two-step approaches.

    \subsection{Convergence analysis of the adaptive  proximal weighted gradient stochastic descent algorithm}
    
        In this section, we demonstrate through simulation the convergence properties of the proposed adaptive  proximal weighted gradient stochastic descent algorithm to optimize for each value of the regularization parameter $\lambda$ in a grid $\Lambda$    the penalized log-likelihood criterion defined in Equation  (\ref{eq:penalized_estimator}). To that purpose we consider a simple toy model where the value of the parameter $\theta$ that maximizes  the targeted criterion 
    can be calculated exactly though       an explicit expression   in order to  highlight the accuracy of the proposed algorithm. Let us therefore consider  for $\nobs \in  \mathbb{N}^*$ and $J \in  \mathbb{N}, J \geq 2$ the following linear mixed effects model where for all $i \in \{1, ..., \nobs\}$ 
       
            \begin{equation}\label{eq:lmem_calcul_exact_vect} 
                \stackeq[rl]{
                    \obs_{i} &= \covgd_{i}\beta +  W \lat_{i}  +\varepsilon_{i}\\
                    \lat_{i}  &\simiid\mc N(0,\gamma^2 )\\
                    \varepsilon_{i} &\simiid\mathcal{N}(0,\sigma^2 I_J)
                }
            \end{equation}
            where $\obs_{i} =(\obs_{i,1}, ..., \obs_{i,J}) \in \mathbb{R}^J$, $\covgd_{i}=(\covgd_{i,1}, ..., \covgd_{i,J}) \in \mathcal{M}_{J,p}(\mathbb{R})$, $\varepsilon_{i} =(\varepsilon_{i,1}, ..., \varepsilon_{i,J}) \in \mathbb{R}^J$ and $W $ is the vector of $\mathbb{R}^J$ with all components equal to $1$.  We assume  that for each $i$ the vectors $\lat_{i} $ and $\varepsilon_{i} $ are independent. Thus the marginal distribution of $\obs_{i}$ is Gaussian with expectation $\covgd_{i}\beta$ and covariance matrix $\Gamma=\sigma^2 I_J+ \gamma^2 W W^t$. 
           We assume also that the random vectors  $(\lat_{i} )$ are independent., as well as  the random vectors $(\varepsilon_{i} )$.
    
%

        Note that thanks to the Woodbury formulae,  we get an explicit expression for the inverse of $\Gamma$, given by:
        \begin{equation}
        \label{eq:invgamma}
            \Gamma^{-1}=\frac{1}{\sigma^2} I_J -  \frac{\gamma^2}{\sigma^2 (  J \gamma^2 +  \sigma^2  ) } W W^t
        \end{equation}
                
         We assume that the variance parameters $\gamma^2$ and $\sigma^2$ are known and  that  parameters  $\beta \in \setr^\dimcovgd$ are estimated. Thus the penalized criterion writes, up to an additive constant:
        
         \begin{eqnarray}
             \mathcal{C}_\lambda(\beta)&=&\log g(\beta ; Y)- \regularization \som k1p \abs{\beta_k} \\
             &=&- \frac{1}{2} \sum_{i=1}^\nobs  (\obs_{i} -\covgd_{i}\beta)^t  \Gamma^{-1}  (\obs_{i} -\covgd_{i}\beta) - \regularization \som k1p \abs{\beta_k} 
        \end{eqnarray}
        
             Using the expression given in  (\ref{eq:invgamma}), we get:
            \begin{eqnarray*}
             \mathcal{C}_\lambda(\beta) &=&- \frac{1}{2\sigma^2} \sum_{i=1}^\nobs  \|\obs_{i} -\covgd_{i}\beta \|^2  - \frac{\gamma^2}{2\sigma^2 (  J \gamma^2 +  \sigma^2  ) } \sum_{i=1}^\nobs 
             \left(  \sum_{j=1}^J   \obs_{ij} - (\sum_{j=1}^J \covgd_{ij})\beta \right)^2 
           - \regularization \som k1p \abs{\beta_k} 
         \end{eqnarray*}
         where $X_{ij}$ is the $j$-th line of matrix $X_i$.
         Let us define the assumption $(H1)$: $\forall 1 \leq  i \leq \nobs$ and $\forall 1 \leq l \leq p$ $\sum_{j=1}^J \covgd_{ijl}=0$. If the design matrix $X$ satisfies $(H1)$, then the associated penalized criterion writes, up to an additive constant:
         \begin{eqnarray*}
             \mathcal{C}_\lambda(\beta) &=&- \frac{1}{2\sigma^2} \sum_{i=1}^\nobs  \|\obs_{i} -\covgd_{i}\beta \|^2  
           - \regularization \som k1p \abs{\beta_k} 
         \end{eqnarray*}
         Note that this expression of the criterion corresponds exactly  to the  LASSO criterion expression in the classical homoscedastic linear regression. If we assume additionally that the design $X$ is orthogonal, meaning that  $X^tX=I_p$, we get an explicit expression of the value $  \hat\beta^{PEN}_\regularization$ of $\beta$ that maximizes this criterion, depending on the expression of the ordinary least square estimate denoted by $   \hat\beta^{OLS}$ (see \cite{tibshirani_regression_1996}): 
         \begin{equation}\label{eq:penalized_estimator_beta_only_expression}
                   \hat\beta^{PEN}_\regularization= sgn(\hat{\beta}^{OLS})(  |  \hat{\beta}^{OLS}) - \regularization  |  )_+,
                \end{equation}
                  with $\textit{sgn}(x)$ equals $1$ if $x>0$, $-1$ if $x<0$ and $0$ if $x=0$ and      where the penalized estimate is defined as 
        \begin{equation}\label{eq:penalized_estimator_beta_only}
                  \hat\beta^{PEN}_\regularization= \argmax_{\beta} \left\{\log \marg(\beta;\obs) - \pen_\regularization(\beta)\right\} 
                \end{equation}
             and where 
        \begin{equation}\label{eq:OLS}
                    \hat\beta^{OLS}= (X^t  X)^{-1} X^t  Y = X^t  Y  .
                \end{equation}
                
        We detail here how to build in practice an orthogonal  design matrix $X$ satisfying $(H1)$. First we sample for $1 \leq  i\leq \nobs$ a matrix $U_i$ of size $J \times p$ with independent standard Gaussian components. We then centered these matrices by column and define the matrices $(V_i)$ of size $J \times p$ by $V_i=U_i - W m_i$ where $m_i $ is the vector of size $p$ such that $m_{il}= \sum_{j=1}^J U_{ijl}/J $ for each $1 \leq l \leq p$. We denote by $T$ the matrix of size $p \times p$ defined by $T=V^t V$ and consider its Cholesky decomposition $T=L L^t$ where $L$ is a lower triangular matrix. Finally we define the design $X$ by $X=V(L^t)^{-1}$. It is easy to check that this design is orthogonal by construction and satisfies $(H1)$ by construction.

\begin{figure}[!ht]
            \centerline{
            \includegraphics[width = 0.9\textwidth]{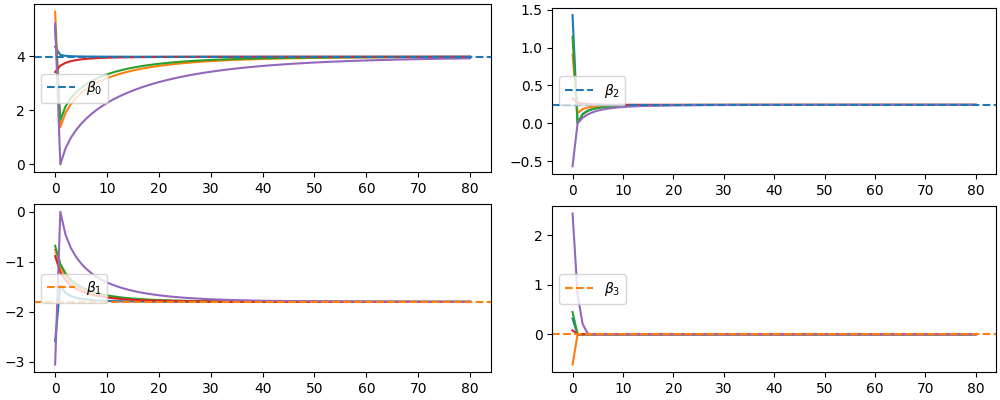}}
             \caption{Trajectories of $5$ runs of the algorithm performed  on the same dataset for the four first components of $\beta$ with randomly chosen initializations and  regularization parameter $\lambda$ equal to $0.6$ leading to over selection of the support.}
             \label{fig:PEN_5runs_over}
         \end{figure}

\begin{figure}[!ht]
            \centerline{
             \includegraphics[width = 0.9\textwidth]{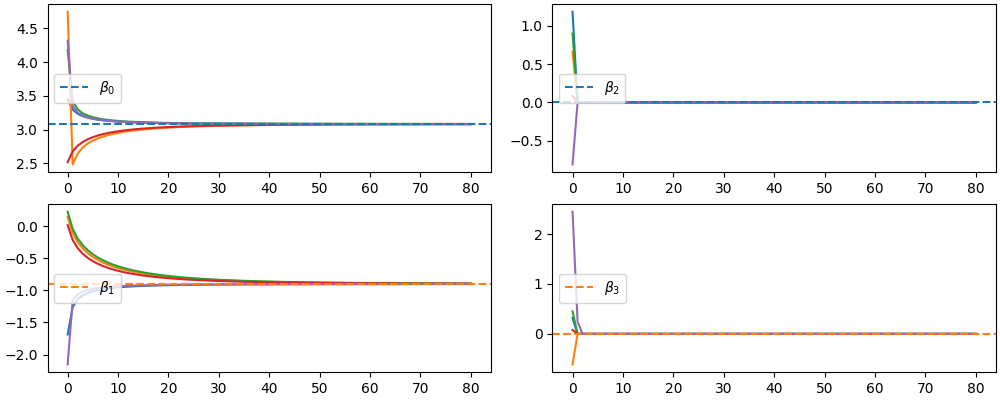}}
            \caption{Trajectories of $5$ runs of the algorithm performed  on the same dataset for the four first components of $\beta$ with randomly chosen initializations and  regularization parameter $\lambda$ equal to $1$ leading to exact selection of the support.}
             \label{fig:PEN_5runs_exact}
         \end{figure}

\begin{figure}[!ht]
            \centerline{
             \includegraphics[width = 0.9\textwidth]{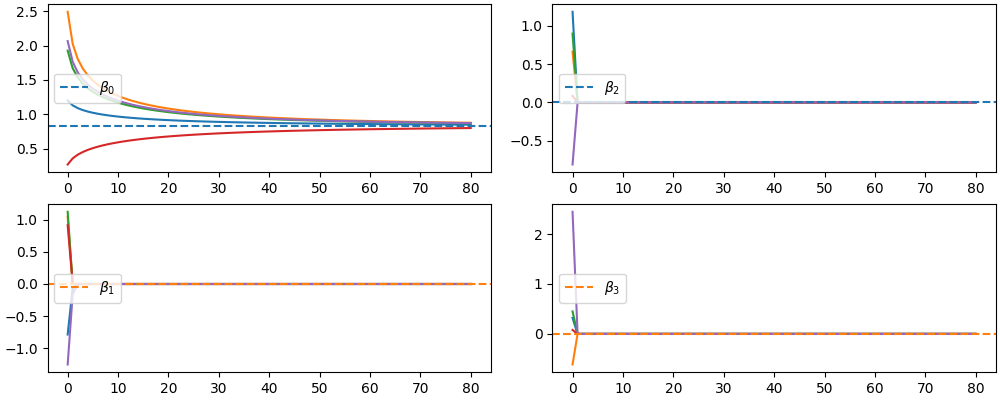}
             }
             \caption{Trajectories of $5$ runs of the algorithm performed  on the same dataset for the four first components of $\beta$ with randomly chosen initializations and  regularization parameter $\lambda$ equal to $2$ leading to under selection of the support.}
             \label{fig:PEN_5runs_under}
         \end{figure}

        Given this design $X$, we generate one dataset  according to the model presented in \cref{eq:lmem_calcul_exact_vect} with $\nobs=100$, $J=5$, $\gamma^2=4^2$, $\sigma^2=2^2$, $p=200$. We choose for the components of $\beta$ the values $4$ and $-3$ for the first and second components respectively, whereas the others are fixed equal to zero. 
        We run the proposed algorithm for three different values of the regularization parameters $\regularization$, namely one leading to select more variables than in the true model, one leading to select the same variables as in the true model, and one leading to select less variables than  in the true model. For each configuration, we run  the algorithm several times with different initializations chosen randomly. We present  five trajectories of the estimates obtained with the  algorithm for the fourth first components of $\beta$ to highlight the convergence toward the target in different settings of regularization: settings with more selected variables   than those in the true support corresponding to over selection case on Figure \ref{fig:PEN_5runs_over}, settings with the same  selected variables   than the  variables in the true support corresponding to exact selection case on Figure \ref{fig:PEN_5runs_exact}, settings with less selected variables  than the variables in the true support corresponding to under selection case on Figure \ref{fig:PEN_5runs_under}. We emphasize again that in this example the true values of the penalized estimates for $\beta$ are calculated explicitely though Equation (\ref{eq:penalized_estimator_beta_only_expression}).  These exact values are represented in dotted horizontal  lines on the graphs. We observe in all settings that the algorithm's trajectories converge with high accuracy  toward the targeted values, which differ depending on the regularization parameter value.

%
%
%

 \subsection{Variable selection in high dimension and re-estimation analysis in a linear mixed effects model}

    For this section and the next one, we generate $n_{\text{run}}$ independent data sets, and we fit the corresponding model using the routine described in \cref{algo:metho} and get an estimate $\hat\theta^{\text{MLE}}_{i}$ for each run $i \in\{1,..., n_{\text{run}}\}$. To compare the results across different scenarios, we calculate several metrics.   We  evaluate  the sensitivity, specificity, and accuracy to study the selection capacity of the method. Sensitivity (Se) measures the proportion of true positives (TP) correctly identified, while Specificity (Sp) quantifies the proportion of true negatives (TN) correctly identified. Accuracy (Ac) represents the overall proportion of correctly classified instances, including both TP and TN.  We abbreviated false negatives and false positives  respectively by FN and FP.

    \begin{equation*}
        \text{Se} = \frac{TP}{TP + FN} ~~
        \text{Sp} = \frac{TN}{TN + FP} ~~
        \text{Ac} = \frac{TP + TN}{P+N}
    \end{equation*}
    
    We also evaluate     the Relative Root Mean Square Errors (RRMSE) to measure the estimation quality of the method defined as:
    \begin{equation*}\label{eq:rrmse_se_sp_ac}
        \text{RRMSE}(\hat\theta^{\text{MLE}}_{k}) = \sqrt{\frac 1{n_{\text{run}}}\sum_{i=1}^{n_{\text{run}}} \frac{(\hat\theta^{\text{MLE}}_{i,k} - \theta_k)^2}{\theta_k^2}} ~~
        \text{mse}(\hat\theta^{\text{MLE}}_{k}) = \frac 1{n_{\text{run}}}\sum_{i=1}^{n_{\text{run}}} (\hat\theta^{\text{MLE}}_{i,k} - \theta_k)^2
    \end{equation*}
     where $\hat\theta^{\text{MLE}}_{i,k}$ is the $k$-th component of $\hat\theta^{\text{MLE}}_{i}$.
    \renewcommand{\varlat}{\gamma} \renewcommand{\vareps}{\sigma}
        
        In this section, we compare our method  \NLMEMLASSO  to the R package \glmmLasso(\cite{GLMMLasso_groll_2014}). We  use the same procedure to select the regularization parameter. We consider the following linear model :

        \begin{equation}\label{eq:lmem_simu_fct} 
            \stackeq[rl]{
                \obs_{i,j} &= \lat_{i,j1}  + \lat_{i2}  \covrep_{j} +\varepsilon_{i,j}\\
                \lat_{i,j,1}  &\siminde\mc N(\mu_1+\covgd_{i,j}\beta,\varlat_1^2)\\
                \lat_{i,2}  &\simiid\mc N(\mu_2,\varlat_2^2)\\
                \varepsilon_{i,j} &\simiid\mathcal{N}(0,\vareps^2)
            }
        \end{equation}
        The model parameters are $\beta \in \setr^\dimcovgd$, $\mu = \p{\mu_1,\mu_2} \in\setr^2$, and $\varlat_1^2,\varlat_2^2, \vareps^2 \in \p{\setr^*_+}^3$. 
        In this model, high-dimensional covariates model variability between individuals. 
        We generated $100$ data sets independently according to \cref{eq:lmem_simu_fct} for several scenarios. It is assumed that each individual is observed $J=10$ times, at the same instants spread uniformly over a range between $0$ and $1$. We use the following values for the parameters : $\mu = \p{2, 5}$, $\varlat_1 = 1, \varlat_2 = 2$ and $\vareps = 1$. For each different value of $\dimcovgd$, we choose the vector $\beta$ such that the first three components are equal to $(\beta_1, \beta_2,\beta_3) = (8, -10, 20)$ and the rest are equal to zero. Additionally, we generate the matrix of covariates $\covgd$ following a uniform distribution $\covgd_{ijk} \sim \mc U([-1,1]) ; ~\forall i \in \{1, ..., \nobs\}, \forall j \in \{1, ..., \nrep\}, \forall k \in \{ 1, ..., \dimcovgd\}$.
 
        By testing scenarios with increasing numbers of individuals $N\in\{100,200\}$, we want to assess how both methods 
         estimate more or less accurately when the sample size increases. In parallel we show different scenarios where the number of covariates increases $p\in\{200,1000\}$. We assess how both methods, despite the high-dimensional context, manage to select the most explanatory variables. To evaluate the selection variable capacity of both  methods, we present  selection scores in Table \ref{chap3:tab:SeSpAc_LMEM} and proportions of well-selected variables regarding the support on Figure \ref{chap3:fig:support_LMEM}. Whatever the configuration, the proposed procedure \NLMEMLASSO and the procedure \glmmLasso perform very similar regarding the considered indicators of average sensitivity, specificity and accuracy. Let's take a closer look at the performance of the two methods in terms of support selection: both generally find the right support.  However, the procedure \NLMEMLASSO   seems to over-select support less often than the procedure \glmmLasso. Tables  \ref{chap3:tab:LMEM_RMSE} and \ref{tab:LMEM_RMSE_P500} (in Appendix) show the mean value over the $100$ data sets of the Relative Root Mean Square Errors (RRMSE) for each model parameter, along with their estimates obtained after variable selection in the reduced model.    Regarding the parameter's re-estimation in the reduced model,  we emphasize that  both methods give similar good results regarding estimation of regression parameters and of the residual variance; nevertheless  the procedure \NLMEMLASSO also adequately estimates the variances of the random effects, unlike the procedure \glmmLasso. This might be explained by  the tendency of the \glmmLasso method to over-select the support.

        \begin{tableht}
            \caption{Average of sensitivity (Se), specificity (Sp), accuracy (Ac), mean square error (mse) over $100$ repetitions within the LMEM using the \NLMEMLASSO and \glmmLasso methods.}


  \newcommand{\algo}{\scriptsize \NLMEMLASSO & \scriptsize \glmmLasso}
  \resizebox{0.8\textwidth}{!}{
  \begin{tabular}{@{}l*{6}{cc}@{}} 
    \toprule
    & \multicolumn{4}{c}{\bf N = 100} 
    & \multicolumn{4}{c}{\bf N = 200} \\
     \cmidrule(r){2-5} \cmidrule(r){6-9}
    & \multicolumn{2}{c}{ \bf P = 200 }  & \multicolumn{2}{c}{ \bf P = 1000 }
    & \multicolumn{2}{c}{ \bf P = 200 } & \multicolumn{2}{c}{ \bf P = 1000 }\\
    \cmidrule(r){2-3} \cmidrule(r){4-5}
    \cmidrule(r){6-7} \cmidrule(r){8-9}
    & \algo & \algo &  \algo & \algo  \\ 
      \midrule    
  Ac & 1.000 & 1.000 & 1.000 & 1.000 & 1.000 &  1.000 & 1.000 & 1.000 \\ 
  Se & 1.000 & 1.000 & 1.000 & 1.000 & 1.000 & 1.000 & 1.000 & 1.000 \\ 
  Sp & 1.000 & 1.000 & 1.000 & 1.000 &  1.000 & 1.000 & 1.000 & 1.000 \\ 
  mse & 0.005 & 0.005  & 0.001 & 0.001 & 0.003 & 0.002 &  0.001 & 0.001 \\ 
    \bottomrule
  \end{tabular}}
 \label{chap3:tab:SeSpAc_LMEM}
        \end{tableht}

        \begin{figureht}
            \includegraphics[width = 0.5\textwidth]{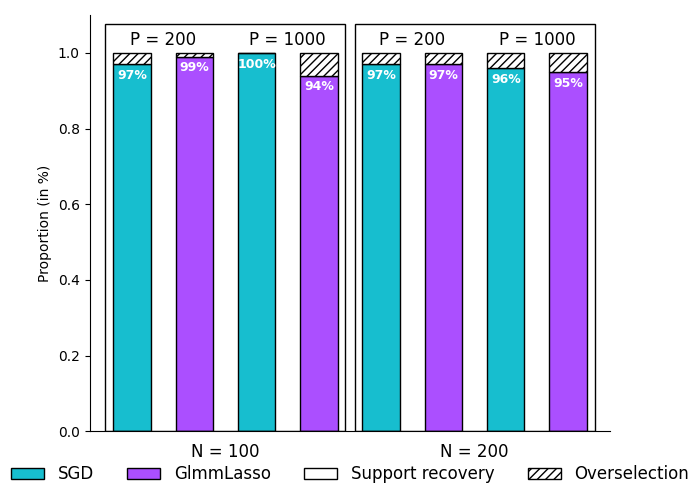}
            \caption{Proportion of correct support selection (solid bars) and support overselection (striped bars) using both procedures \NLMEMLASSO and \glmmLasso for $N \in \{100,200\}$ and for $P \in \{200,  1000\}$.}
              \label{chap3:fig:support_LMEM}
        \end{figureht}

        \begin{tableht}
            \caption{Average Relative Root Mean Square Errors (RRMSE) and parameter estimates over $100$ repetitions for the Linear Mixed Effects Model (LMEM) where $\hat\theta_\NLMEMLASSO$ and $\hat\theta_\glmmLasso$ are the estimates in the reduced model after variable selection  using the \NLMEMLASSO and \glmmLasso methods, respectively.}


    \newcommand{\rrmse}{\scriptsize RRMSE}
    \resizebox{0.8\textwidth}{!}{
    \begin{tabular}{@{}ccc*{4}{cc}@{}} 
        \toprule
           & & & \multicolumn{4}{c}{ \bf P = 200 } & \multicolumn{4}{c}{ \bf P = 1000 } \\  
            \cmidrule(r){4-7} \cmidrule(r){8-11}
           & & & \multicolumn{2}{c}{ \bf \NLMEMLASSO } & \multicolumn{2}{c}{ \bf \glmmLasso }
            & \multicolumn{2}{c}{ \bf \NLMEMLASSO } & \multicolumn{2}{c}{ \bf \glmmLasso } \\
            \cmidrule(r){4-5} \cmidrule(r){6-7} \cmidrule(r){8-9} \cmidrule(r){10-11}
          &  & $\theta^*$ & 
            $\hat\theta$ & \rrmse & $\hat\theta$ & \rrmse & 
            $\hat\theta$ & \rrmse & $\hat\theta$ & \rrmse \\
            \midrule
        \multirow{8}{*}{ \bf N = 100 } 
 & $\mu_1$ & 2.00 & 1.99 & 3.15 & 1.98 & 6.13 & 1.99 & 3.07 & 1.98 & 6.14 \\ 
   & $\mu_2$ & 5.00 & 5.02 & 1.92 & 5.03 & 3.70 & 5.02 & 1.91 & 5.03 & 3.68 \\ 
   & $\gamma_1^2$ & 1.00 & 0.99 & 19.96 & 0.61 & 39.74 & 0.98 & 19.28 & 0.61 & 39.40 \\ 
   & $\gamma_2^2$ & 4.00 & 3.93 & 17.10 & 0.88 & 78.10 & 3.90 & 17.01 & 0.88 & 78.12 \\ 
   & $\sigma^2$ & 1.00 & 1.00 & 5.28 & 0.98 & 5.07 & 1.00 & 5.29 & 0.98 & 5.15 \\ 
   & $\beta_{1}$ & 8.00 & 8.03 & 7.20 & 8.05 & 7.20 & 8.05 & 7.23 & 8.06 & 7.74 \\ 
   & $\beta_{2}$ & -10.00 & -9.95 & 5.49 & -9.97 & 5.81 & -9.87 & 6.04 & -9.92 & 6.30 \\ 
   & $\beta_{3}$ & 20.00 & 19.99 & 2.74 & 19.98 & 2.89 & 19.91 & 2.65 & 19.86 & 2.80 \\ 
            \midrule
        \multirow{8}{*}{ \bf N = 200 }
 & $\mu_1$ & 2.00 & 2.00 & 1.95 & 1.99 & 3.65 & 2.00 & 1.83 & 2.00 & 3.65 \\ 
   & $\mu_2$ & 5.00 & 5.00 & 1.85 & 5.01 & 3.49 & 5.00 & 1.84 & 5.01 & 3.49 \\ 
   & $\gamma_1^2$ & 1.00 & 0.99 & 14.63 & 0.61 & 39.12 & 0.98 & 14.09 & 0.61 & 38.97 \\ 
   & $\gamma_2^2$ & 4.00 & 4.06 & 12.77 & 0.89 & 77.83 & 4.04 & 12.41 & 0.89 & 77.80 \\ 
   & $\sigma^2$ & 1.00 & 1.00 & 4.00 & 0.98 & 4.07 & 1.00 & 4.09 & 0.98 & 3.94 \\ 
   & $\beta_{1}$ & 8.00 & 8.03 & 4.52 & 8.02 & 4.59 & 7.91 & 5.01 & 7.92 & 5.21 \\ 
   & $\beta_{2}$ & -10.00 & -9.94 & 3.70 & -9.96 & 3.85 & -9.88 & 4.12 & -9.90 & 4.11 \\ 
   & $\beta_{3}$ & 20.00 & 20.01 & 1.93 & 20.01 & 1.92 & 20.02 & 2.13 & 20.02 & 2.13 \\ 
        \bottomrule
    \end{tabular}}
 \label{chap3:tab:LMEM_RMSE}
        \end{tableht}


   \subsection{Variable selection in high-dimension and re-estimation analysis in a logistic non-linear mixed-effects model}
                
        We study in this section the model presented in the first example in Equation  (\ref{eq:mem_orangetrees}), where  
        we explain a part of the variability of the logistic’s midpoint by the high-dimensional covariates.
        The model can be written as follows:        
            \begin{equation}\label{eq:nlmemsimu_fct}
                \stackeq[rl]{
                    \obs_{i,j}&= \dfrac{\lat_{i1}}{1 + \exp\p{-\frac{\covrep_{ij} - \lat_{i2}}{\alpha}}}+\varepsilon_{i,j}\\
                    \lat_{i,1} & \simiid \mc N(\mu_1,\varlat_1^2)\\
                    \lat_{i,2} &\siminde \mc N(\mu_2 +\covgd_i\beta,\varlat_2^2)\\
                    \varepsilon_{i,j}&\simiid\mathcal{N}(0,\vareps^2)
                }
            \end{equation}        
        The model parameters are $\alpha \in \setr,  \beta \in \setr^\dimcovgd$, $\mu = \p{\mu_1,\mu_2} \in\setr^2$, and $\varlat_1^2,\varlat_2^2, \vareps^2 \in \p{\setr^*_+}^3$. The quantity $\lat_{i1}$ represents the asymptotic maximum value of the curve, $\lat_{i2}$ represents the value of the logistic’s midpoint, and $\alpha$ represents the logistic growth rate.  The two first parameters are modeled at the individual level of the mixed model, whereas the third one is modeled at the population level.
        Note that due to the presence of the fixed effect $\alpha$, the joint density of $(Y_i,\lat_i)$ does not belong to the curved exponential family as defined in \cite{delyon1999convergence}. 
        We generate $100$ data sets independently according to \cref{eq:nlmemsimu_fct}  for several scenarios. It is assumed that each individual was observed $J= 15$ times,        
        at the same instants equally spaced over a range between $150$ and $3000$. We use the following values for the parameters: $\mu = \p{200, 1200}$, $\varlat_1^2 =7^2$, $\varlat_2^2= 30^2$, $\alpha = 300$ and $\vareps^2  = 30$. For each different value of $\dimcovgd$, we choose the vector $\beta$ such that the first three components are equal to $(100, 200, -300)$ and the rest equal to zero. Furthermore, we generate the matrix of covariates $\covgd$ with $\nobs$ rows and $\dimcovgd$ columns, following a uniform distribution $\covgd_{i,k} \sim \mc U([-1,1]) ; ~\forall i \in \{ 1, ..., \nobs\}, \forall k \in \{ 1, ..., \dimcovgd\}$.

        As in the previous analysis, we present scenarios with an increasing number of individuals $\nobs \in \{100,200\}$ and different covariates size $p \in \{200,500,1000\}$.     We first  assess the selection variable capacity of our  method, by presenting  selection scores in Table \ref{chap3:tab:SeSpAc_NLMEM} and proportions of well-selected variables regarding the support on Figure \ref{chap3:fig:prop_supp_NLMEM}. We also indicate  over-selection cases. The later shows that as the number of variables $p$ increases, support selection performance deteriorates; in particular, a higher proportion of support over-selection is observed. However, for a fixed number of variables $p$, the support selection performance improves as the number of individuals N increases. 

        Table \ref{chap3:tab:NLMEM_RMSE} summarizes the results of the estimation of the model parameters in the reduced model after variable selection step,  for $N = 100$ and $N = 200$. 
      We observe that the estimates of the parameters are accurate in all configurations. As expected, the relative root mean square errors decrease with $\nobs$ for a given covariate size $p$. For $\nobs=100$, we observe that the relative root mean square errors generally increase with $p$, which can be explained by the fact that increasing the size of the covariates increases the difficulty of the selection procedure. The phenomenon is less pronounced as the sample size increases, reflecting the asymptotically expected behavior.

        \begin{tableht}\centering
            \caption{Average sensitivity (Se), specificity (Sp), accuracy (Ac), and mean square error (mse) using the procedure \NLMEMLASSO over 100 repetitions for the logistic non-linear mixed effects model}


  \begin{tabular}{@{}l*{6}{c}@{}}
    \toprule
    & \multicolumn{3}{c}{ \bf N = 100 } & \multicolumn{3}{c}{ \bf N = 200 } \\
    \cmidrule(r){2-4} \cmidrule(r){5-7}
    & \bf P = 200 & \bf P = 500  & \bf P = 1000  
    & \bf P = 200 & \bf P = 500  & \bf P = 1000  \\ 
      \midrule
        Ac & 0.997 & 0.999 & 0.999 & 1.000 & 1.000 & 1.000 \\ 
        Se & 0.980 & 0.927 & 0.933 & 1.000 & 1.000 & 1.000 \\ 
        Sp & 0.997 & 0.999 & 1.000 & 1.000 & 1.000 & 1.000 \\ 
        mse & 0.819 & 0.851 & 0.387 & 0.127 & 0.057 & 0.026 \\ 
    \bottomrule
    \end{tabular}
 \label{chap3:tab:SeSpAc_NLMEM} 
        \end{tableht}

        \begin{figureht}
            \includegraphics[width = 0.6\textwidth]{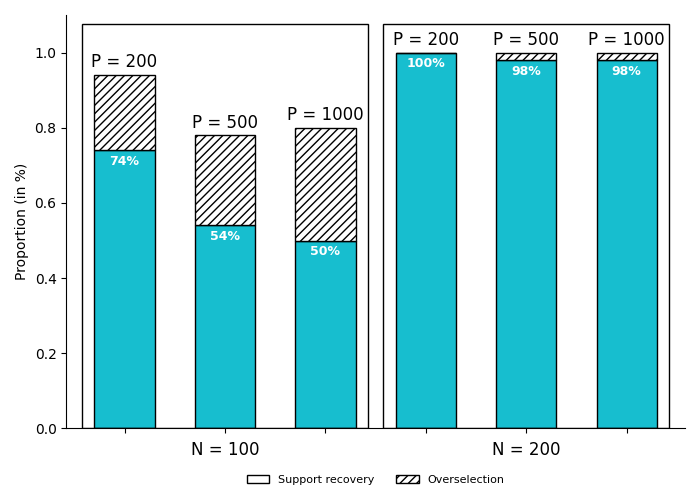}
            \caption{ Proportion of correct support selection (solid bars) and overselection (striped bars) in the logistic non-linear mixed effects model  using the procedure \NLMEMLASSO for $N \in \{100,200\}$ and for $P \in \{200, 500, 1000\}$.}
             \label{chap3:fig:prop_supp_NLMEM}
        \end{figureht}

\begin{tableht}
            \caption{Average Relative Root Mean Square Errors (RRMSE) and parameter estimates over $100$ repetitions for the logistic non-linear mixed effects model  where $\hat\theta$ is the estimate in the reduced model after variable selection  using the \NLMEMLASSO.}


\newcommand{\rrmse}{\scriptsize RRMSE}
    \begin{tabular}{@{}ccc*{3}{cc}c@{}} 
        \toprule
           & & & \multicolumn{2}{c}{ \bf P = 200 } & \multicolumn{2}{c}{ \bf P = 500 } & \multicolumn{2}{c}{ \bf P = 1000 }\\
            \cmidrule(r){4-5} \cmidrule(r){6-7} \cmidrule(r){8-9}
           & & $\theta^*$ & 
            $\hat\theta$ & \rrmse & $\hat\theta$ & \rrmse & $\hat\theta$ & \rrmse & \\ 
            \midrule
        \multirow{9}{*}{ \bf N = 100 } 
            & $\mu_1$ & 200.00 & 199.90 & 0.38 & 199.90 & 0.40 & 199.93 & 0.39 \\ 
            & $\mu_2$ & 1200.00 & 1200.02 & 0.35 & 1199.87 & 0.38 & 1200.20 & 0.37 \\ 
            & $\gamma_1^2$ & 49.00 & 48.46 & 17.97 & 48.20 & 17.68 & 48.31 & 17.82 \\ 
            & $\gamma_2^2$ & 900.00 & 861.29 & 38.85 & 1062.15 & 63.75 & 1057.36 & 67.47 \\ 
            & $\alpha$ & 300.00 & 299.86 & 0.77 & 299.99 & 0.71 & 300.16 & 0.76 \\ 
            & $\sigma^2$ & 30.00 & 30.06 & 5.74 & 30.09 & 5.82 & 30.09 & 5.77 \\ 
            & $\beta_{0}$ & 120.00 & 119.90 & 3.47 & 119.38 & 3.85 & 118.44 & 3.92 \\ 
            & $\beta_{1}$ & 70.00 & 69.74 & 5.86 & 69.27 & 6.80 & 69.68 & 6.75 \\ 
            & $\beta_{2}$ & 40.00 & 37.81 & 26.31 & 31.32 & 47.85 & 31.56 & 45.43 \\ 
            \midrule
        \multirow{9}{*}{ \bf N = 200 } 
            & $\mu_1$ & 200.00 & 199.97 & 0.26 & 199.99 & 0.25 & 199.99 & 0.26 \\ 
            & $\mu_2$ & 1200.00 & 1199.75 & 0.23 & 1199.94 & 0.22 & 1199.96 & 0.23 \\ 
            & $\gamma_1^2$ & 49.00 & 48.76 & 11.83 & 48.74 & 12.06 & 48.71 & 11.95 \\ 
            & $\gamma_2^2$ & 900.00 & 894.39 & 14.65 & 891.41 & 14.39 & 890.88 & 14.16 \\ 
            & $\alpha$ & 300.00 & 299.82 & 0.50 & 299.91 & 0.50 & 299.91 & 0.51 \\ 
            & $\sigma^2$ & 30.00 & 30.10 & 3.76 & 30.10 & 3.97 & 30.12 & 3.86 \\ 
            & $\beta_{0}$ & 120.00 & 119.87 & 2.35 & 119.88 & 2.53 & 119.35 & 2.42 \\ 
            & $\beta_{1}$ & 70.00 & 69.97 & 4.50 & 70.37 & 4.50 & 69.13 & 4.36 \\ 
            & $\beta_{2}$ & 40.00 & 40.53 & 6.89 & 39.39 & 7.11 & 40.29 & 6.59 \\
        \bottomrule
    \end{tabular}
 \label{chap3:tab:NLMEM_RMSE}
        \end{tableht}

\begin{figure}[!ht]
            \centerline{\includegraphics[width = 0.75\textwidth]{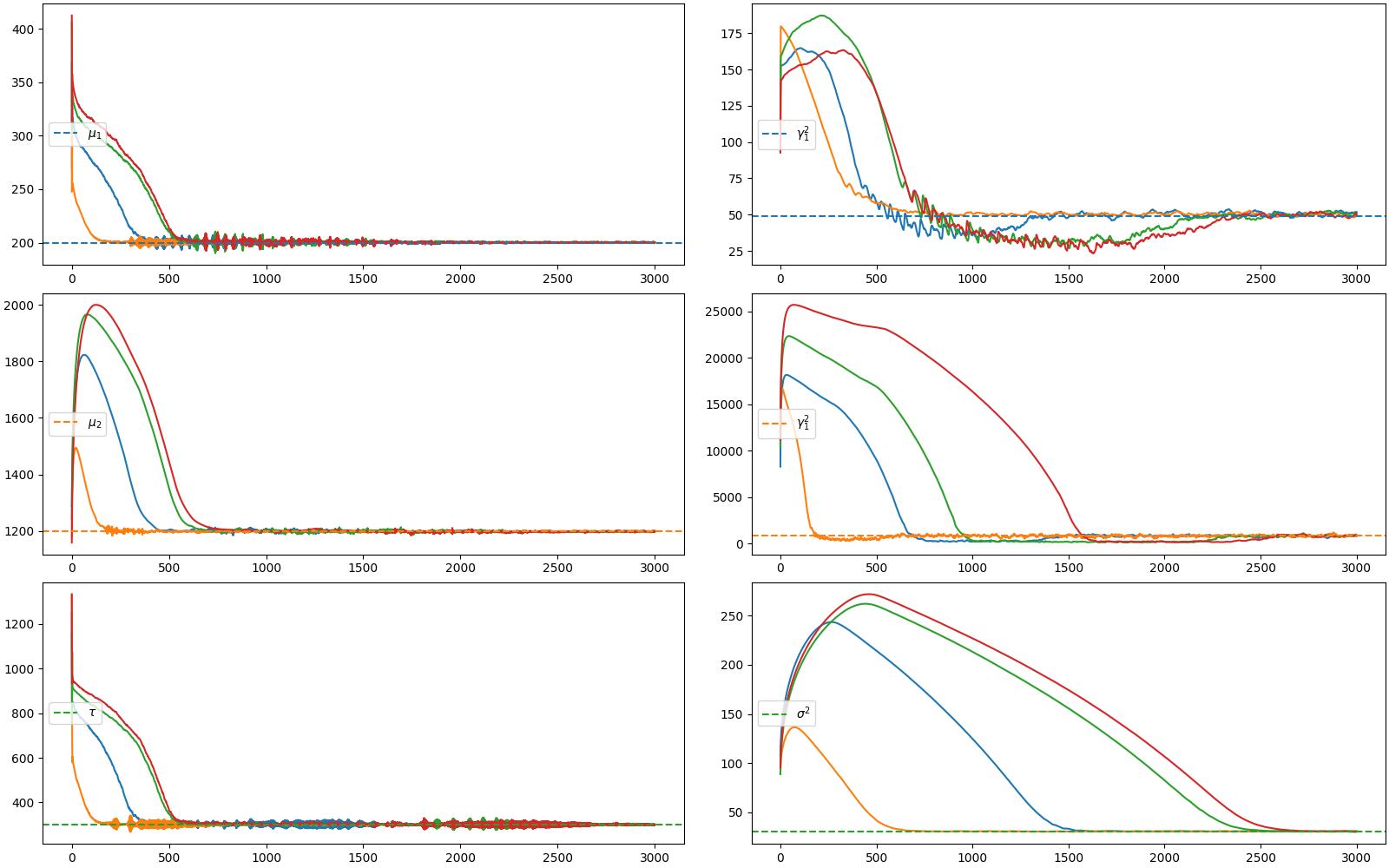}}
            \caption{Parameter estimates in the reduced model across the iterations of the adaptive stochastic gradient algorithm within the logistic mixed model.  Dotted lines : true values}
            \label{fig:PSPG_iter}
        \end{figure}

        As an illustration of the behavior of the optimization algorithm in the reduced model, Figure \ref{fig:PSPG_iter} displays the estimated parameter as a function of iterations during the execution of the adaptive stochastic gradient algorithm. We can observe the good convergence of the algorithm toward a limit which is very closed to the true parameter value, as expected for the maximum likelihood estimate.
        

    \subsection{Comparison of integrated and two steps approaches for variable selection in a non-linear mixed effects model}

        In this section,  following the study performed in  \citep{naveau2024bayesian}, we  compare the performance for variable selection of three procedures. The first procedure  consists in two steps,  namely first in treating the observations of each individual separately by ajusting one non linear regression by individual to get estimates of its parameters and second performing a high dimensional linear regression  between these estimates and the covariates to select the more relevant. The second procedure consists also in two steps, namely first adjusting a non linear mixed model to estimate the individual parameters of all individuals simultaneaously, and second performing a high dimensional linear regression  between these estimates and the covariates to select the more relevant. The third procedure is the one proposed in  section  \ref{algo:metho} and integrates both steps in one single step.

    We  carry out the first procedure called \nlm-\glmnet  using the method of least squares implemented in the \nlm \R function \citep{Dennis1996NLOptimization} to estimate   individual parameters and then  perform variable selection with the LASSO procedure using the \glmnet package \citep{Friedman2010glmnet}. This last step is performed using  cross-validation in order to select an optimal regularization parameter chosen such that the mean squared error  is within one standard deviation of the minimum mse as recommended by \citet{hastie2009statistical}.
     The second procedure called \saemix-\glmnet uses the 
     Stochactic Approximation Expectation Maximization (SAEM) algorithm \citep{Kuhn_MCMCSAEM_2004} implemented in the R package \saemix \citep{saemix2017Comets} to adjust the population mixed model and obtain the individual parameters estimates. The variable selection is then performed as in the first procedure. The third procedure called  \NLMEMLASSO is performed as detailed in section  \ref{algo:metho}.

       For this comparison, we use the pharmacodynamic non linear mixed effect model. presented in the second example in Equation (\ref{eq:mem_theophylline}).  We generate 100 independent datasets according to the following model with $\nobs = 200$ individuals and $p = 500$ covariates :
        \begin{equation}\label{chap3:eq:pharmaco_model}
            \stackeq[rllllll]{
            \obs_{ij} &= \dfrac{D \lat_{i1}}{\mc V(\lat_{i1} - \lat_{i2})} 
                \left( \exp\left(-\dfrac{\lat_{i2}}{\mc V} \covrep_{ij}\right) - \exp\left(-\lat_{i1} \covrep_{ij}\right) \right) + \eps_{ij}, \\[2ex]
            \lat_i &= \meanlat +  \covgd_i \beta+ \xi_i, \quad
    \xi_i \simiid \mc N_2\left(0, 
            \begin{pmatrix}
            \varlat_1^2 & 0 \\
            0 & \varlat_2^2
            \end{pmatrix}
            \right), \quad
            \eps_{ij} \simiid \mc N(0, \vareps^2),
            }
        \end{equation}
        The model parameters are $\meanlat = (\meanlat_1, \meanlat_2) \in \setr^2$, $\beta\in\mc M_{p,2}$ and $(\varlat_1^2, \varlat_2^2, \vareps^2) \in \p{\setr^*_+}^3$.   The quantity $\mc V$ is the distribution volume,  whereas $D$ is the administered dose. These two parameters are supposed to be known. The quantity $\lat_{i1}$ represents the absorption rate, $\lat_{i2}$ represents the elimination rate.       We consider  the following observations times $[0.05, 0.15, 0.25, 0.4, 0.5, 0.8, 1, 2, 7, 12, 24, 40]$. These times are chosen to cover a wide range of the pharmacokinetic profile.
        
       We generate datasets with a censored  rate $\rho_{censored}$ varying in  $\{0\%, 10\%, 20\%, 30\%, 40\%\}$.  More precisely a proportion $\rho_{censored}$ of the individuals are observed at only the first three times and supposed to leave the study after that. The other individuals are observed at all times. We use the following values for the parameters: $\meanlat = (6,8)$, $\varlat_1 = 0.2$, $\varlat_2 = 0.1$ and $\vareps = 10^{-3}$. The high-dimensional parameter $\beta$ is chosen such that :
        \begin{equation*}
            \beta = 
            \begin{pmatrix}
                3 & 2 & 1 & 0 & 0 & 0 & \dots & 0\\
                0 & 0 & 3 & 2 & 1 & 0 & \dots & 0
            \end{pmatrix}^T
        \end{equation*}
        The individual covariates $\covgd_i$ are generated from a binomial distribution with a probability of success of $0.2$ and are standardized.    We emphasize here that the second procedure \saemix-\glmnet is carried out on   the following mixed model without high-dimensional covariates presented in Equation (\ref{chap3:eq:pharmaco_model_nogd})  to estimate the individual parameter $\lat_i$ with the package \saemix in the first step. 
        
        \begin{equation}\label{chap3:eq:pharmaco_model_nogd}
            \stackeq[rllllll]{
            \obs_{ij} &=  \dfrac{D \lat_{i1}}{\mc V(\lat_{i1} - \lat_{i2})} 
                \left( \exp\left(-\dfrac{\lat_{i2}}{\mc V} \covrep_{ij}\right) - \exp\left(-\lat_{i1} \covrep_{ij}\right) \right) + \eps_{ij}, \\[2ex]
            \lat_i &= \meanlat + \xi_i,  \quad
            \xi_i \simiid \mc N_2\left(0, 
            \begin{pmatrix}
            \varlat_1^2 & 0 \\
            0 & \varlat_2^2
            \end{pmatrix}
            \right), \quad
            \eps_{ij} \simiid \mc N(0, \vareps^2),
            }
        \end{equation}

        Table \ref{chap3:tab:prediction_latent} shows the mean estimations  errors (mee) in the first step of  \nlm-\glmnet and \saemix-\glmnet  methods over the 100 datasets. As expected, both methods estimate the first parameter $\lat_1$ with a smaller error than the second parameter $\lat_2$. Moreover, both struggle to estimate the second parameter  $\lat_2$ when the censoring rate increases.  However the mixed-effect approach handle the missing data better than the individual by individual  approach as the error of the second parameter $\lat_2$ is smaller for \saemix-\glmnet.

        \begin{tableht}
        \resizebox{0.7\textwidth}{!}{
            \begin{tabular}{lrrrrrr}
                \toprule
                &  & \multicolumn{5}{c}{Censored data rate} \\ \cline{3-7}
                & Method & 0\% & 10\% & 20\% & 30\% & 40\% \\ 
                \midrule
            \multirow{2}{*}{$mee(\lat_1)$} 
                & \nlm-\glmnet &  0.088 & 0.096 & 0.104 & 0.112 & 0.120 \\ 
                & \saemix-\glmnet & 0.102 & 0.107 & 0.111 & 0.116 & 0.120 \\ 
                            \midrule
            \multirow{2}{*}{$mee(\lat_2)$} 
                & \nlm-\glmnet & 0.143 & 0.657 & 1.188 & 1.701 & 2.222 \\ 
                & \saemix-\glmnet & 0.150 & 0.407 & 0.663 & 0.905 & 1.151 \\ 
                \bottomrule
            \end{tabular}
            }
            \caption{Mean estimation errors (mee) of the parameters $\lat_1$ and $\lat_2$ for the two-step approaches \nlm-\glmnet and \saemix-\glmnet.}
            \label{chap3:tab:prediction_latent}
        \end{tableht}

   \begin{figureht}
            \centerline{\includegraphics[width = 0.9\textwidth]{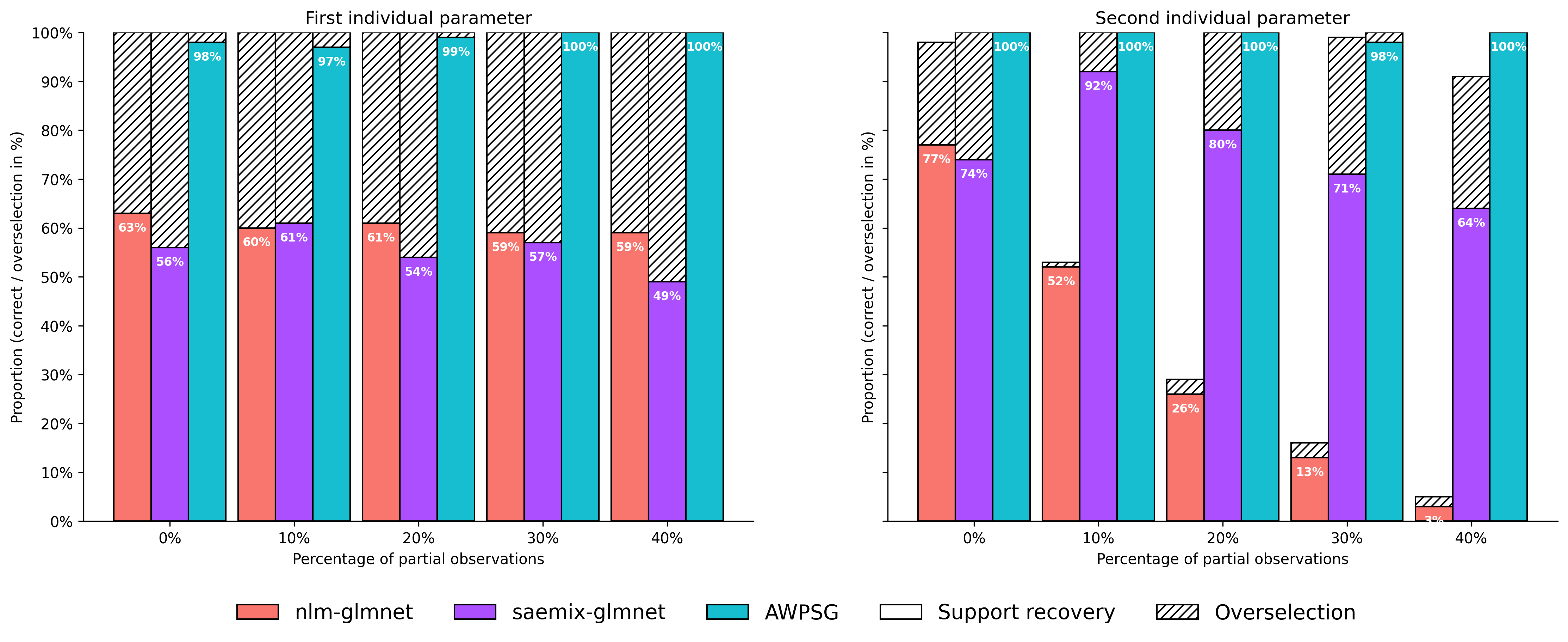}}
            \caption{Proportion of correct support selection (solid bars) and overselection (striped bars) for the \NLMEMLASSO ,  \saemix-\glmnet and \nlm-\glmnet methods under increasing censoring rates.}
            \label{chap3:fig:two_step_comparison}
        \end{figureht}

        Figure \ref{chap3:fig:two_step_comparison} shows the comparison of variable selection performance between the proposed \NLMEMLASSO method and the two-step approaches \nlm-\glmnet and \saemix-\glmnet. On one hand variable selection performs well on the first parameter $\lat_1$ and their associated covariates, and is not impacted by censoring.           On the other hand,   variable selection on  the second parameter is  deteriorated  by censorship in all procedures.  Figure \ref{chap3:fig:two_step_comparison}  highlights the interests of using procedures based on mixed effects models and integrated approach.  The two-step approach \saemix-\glmnet is more successful than the two-step approach \nlm-\glmnet  in finding exact support for different censoring levels.  This can be explained by the lower estimation errors obtained by \saemix-\glmnet observed in Table  \ref{chap3:tab:prediction_latent}.  However, unlike \NLMEMLASSO, it overselects explanatory variables. In all cases, the \NLMEMLASSO method outperforms the two-step approaches.

\section{Conclusion and perspectives}

    In this work, we jointly address variable selection in high-dimension and parameter estimation in general mixed-effect model. We emphasize that the proposed procedure called \NLMEMLASSO can handle well linear and non-linear mixed effects models, without requiring them to belong to the curved exponential family. We perform an adaptive weighted proximal stochastic gradient algorithm to deal simultaneously with the latent variables and the LASSO penalty when maximizing the penalized likelihood. We highlight the well behavior of this optimization algorithm in various settings through simulation. We use the eBIC criterion to perform the choice of the regularization parameter leading to the selection of the variables corresponding to the reduced model. Our integrated methodology for variable selection has been thoroughly evaluated in simulation study to demonstrate its performance in linear and non-linear mixed models. Its performance has also been compared to the ones of the procedure \glmmLasso in the case of linear mixed models. Moreover, we emphasize that the proposed procedure can handle  a wide range of modeling choices for the random effects and the error term, including in particular correlation structure for the random effects. 
    
    There are several interesting perspectives to this work. The first one consists in studying the effect of the presence of correlation between the variables involved in the selection which is known to deteriorate the accuracy of the selection procedure. The second one consists in studying theoretical properties of the proposed procedure, in particular convergence of the  optimization algorithm combining weighted proximal operator and adaptive preconditioned gradient descent and support consistency guarantees.  Finally, it would be of great interest to assort the re-estimation results in the reduced model with confidence region, by
    studying theoretical properties of the estimator and prediction after the variable selection step as done in post model selection inference. 
\section*{Funding and Acknowledgements}

     This work was funded by the https://stat4plant.mathnum.inrae.fr/(Stat4Plant) project ANR-20-CE45-0012.

\bibliographystyle{plainnat}
\bibliography{references}

\newpage \renewcommand\thesubsection{Appendix \Alph{subsection}.}
\section*{Appendix}

    \subsection{Simulation study results for LMEM  }

        \begin{table}[!ht]  \centering
            \caption{Average Relative Root Mean Square Errors (RRMSE) and parameter estimates over $100$ repetitions for the Linear Mixed Effects Model (LMEM)  with $p = 500$ where $\hat\theta_\NLMEMLASSO$ and $\hat\theta_\glmmLasso$ are the estimates in the reduced model after variable selection  using the \NLMEMLASSO and \glmmLasso methods, respectively.}


    \newcommand{\rrmse}{\scriptsize RRMSE}
    \begin{tabular}{@{}cc*{4}{cc}@{}} 
        \toprule             
            & & \multicolumn{8}{c}{\bf P = 500}  \\
            \cmidrule(r){3-10}
            & & \multicolumn{4}{c}{ \bf N = 100 } & \multicolumn{4}{c}{ \bf N = 200 } \\  
            \cmidrule(r){3-6} \cmidrule(r){7-10}
            & & \multicolumn{2}{c}{ \NLMEMLASSO } & \multicolumn{2}{c}{ \glmmLasso }
            & \multicolumn{2}{c}{ \NLMEMLASSO } & \multicolumn{2}{c}{ \glmmLasso } \\
            \cmidrule(r){3-4} \cmidrule(r){5-6} \cmidrule(r){7-8} \cmidrule(r){9-10}
            & $\theta^*$ & 
             $\hat\theta_\NLMEMLASSO$ & \rrmse & $\hat\theta_\glmmLasso$ & \rrmse & 
            $\hat\theta_\NLMEMLASSO$ & \rrmse & $\hat\theta_\glmmLasso$ & \rrmse \\
            \midrule
                $\mu_1$ & 2.00 & 1.98 & 6.36 & 1.98 & 6.15 & 1.98 & 3.97 & 2.00 & 3.64 \\ 
                $\mu_2$ & 5.00 & 5.03 & 4.03 & 5.03 & 3.69 & 5.01 & 2.70 & 5.01 & 3.49 \\ 
                $\gamma_1^2$ & 1.00 & 1.00 & 20.55 & 0.61 & 39.38 & 1.01 & 13.92 & 0.61 & 39.53 \\ 
                $\gamma_2^2$ & 4.00 & 3.93 & 17.32 & 0.88 & 78.11 & 4.03 & 12.72 & 0.88 & 77.92 \\ 
                $\sigma^2$ & 1.00 & 0.99 & 4.99 & 0.98 & 5.00 & 1.01 & 3.91 & 0.98 & 3.97 \\ 
                $\beta_{1}$ & 8.00 & 8.14 & 7.33 & 8.13 & 7.45 & 8.08 & 5.35 & 7.98 & 5.45 \\ 
                $\beta_{2}$ & -10.00 & -9.98 & 6.05 & -10.00 & 6.63 & -9.99 & 3.91 & -9.96 & 4.22 \\ 
                $\beta_{3}$ & 20.00 & 19.98 & 3.30 & 19.97 & 3.42 & 19.98 & 1.91 & 19.97 & 2.14 \\ 
        \bottomrule
    \end{tabular}
 \label{tab:LMEM_RMSE_P500}
        \end{table}

        
\end{document}